\documentclass[11pt]{article}

\usepackage{a4}
\usepackage{latexsym}
\usepackage{amsthm}
\usepackage{amsmath}
\usepackage{amsfonts}
\usepackage{german}

\usepackage{graphics}

\theoremstyle{plain}

\theoremstyle{definition}

\newcommand{\bmx}{\left[\begin{matrix}}
\newcommand{\emx}{\end{matrix}\right]}

\title{A nice involution for multivariable polynomial rings}
\author{Wiland Schmale\thanks{ University of Oldenburg (retired), Germany,  email: wiland.schmale@uni-oldenburg.de}}
\date{25. November 2020}

\begin{document}
\parindent3mm	
\maketitle

\begin{abstract}
The principal minors of the Toeplitz matrix $\left( x_{i-j+1}\right)_{1\le i,j,\le n}$, where $x_0=1, x_k=0$ if $k\le -1$, directly determine an involution of the polynomial ring $R[x_1, ... ,x_n]$ over any commutative ring $R$.
\end{abstract}

{\bf keywords: }  involution, multivariable polynomial ring, principal Toeplitz minors\ \ \ 

\vspace{ .50cm}

In \cite{algcond} principal minors have been used to transform a conjecture on certain Toeplitz pencils into an equivalent algebraic conjecture. A closer look on this transition to minors instead of the original variables shows that it is actually an involution.

Let $R$ be any commutative ring, $R[x_1, ... ,x_n]$ the polynomial ring in the independent variables $x_1, ... ,x_n$, $n\ge 1$,  over $R$ and
\[
T(k)={\small\bmx x_1&1&0&\dots&0\\
          x_2&x_1&1&\dots&0 \\
          \vdots&\vdots && \ddots&\vdots \\
          x_{k-1}&x_{k-2}& \dots& &1\\
          x_{k}&x_{k-1}& \dots& \dots&x_1\\
          \emx,
          m_k=m_k(x_1, ... ,x_k)=\det(T(k))
           \text{ for }1\le k\le n,}
\]

Note that $T(k)$ is a submatrix and $m_k$ a principal minor of $T(n)$. One has $m_1=x_1$, $m_2=x_1^2-x_2$ and so on.

Let now the substitution ring homomorphism $\varphi $ of $R[x_1, ... ,x_n]$ over $R$ be defined as

\[
	\varphi: R[x_1, ... ,x_n] \rightarrow R[x_1, ... ,x_n]	, \ \ x_k \mapsto m_k,\  \text{ for } 1\le k\le n. 
\]

The following can then be observed:

{\bf Involution property: }
\[
	\varphi \circ \varphi =\text{ id}_{R[x_1, ... ,x_n]}= \text{identity map on } R[x_1, ... ,x_n]
\]

The proof will be based on the following recursive properties of principal minors of $T(n)$.

{\bf Recursions for $\pmb{ m_k}$ and $\pmb{ x_k}$: } For $1\le k\le n$

\begin{equation}
	m_k=\sum_{i=1}^k   (-1)^{i-1} x_i m_{k-i},\ \   m_0:=1
\end{equation}
and
\begin{equation} 
	x_k=\sum_{i=1}^k   (-1)^{i-1} m_i x_{k-i},\ \   x_0:=1
\end{equation}

{\bf Proof: } 
In order to obtain (1) we will prove the following slightly more general relation for any first column entries $a_1, ... ,a_k$:

\[
	\det {\small\bmx a_1&1&0&\dots&0\\
          a_2&x_1&1&\dots&0 \\
          \vdots&\vdots && \ddots&\vdots \\
          a_{k-1}&x_{k-2}& \dots& &1\\
          a_{k}&x_{k-1}& \dots& \dots&x_1\\
          \emx}
          =\sum_{i=1}^k  (-1)^{i-1}a_i m_{k-i}
         , \text{ where } \  m_0:=1.
\]
 This implies relation (1) as a special case. The instance ``$k=1$'' beeing trivial we proceed inductively and obtain for $1\le k\le n-1$
 
 \begin{align*}
 	\det {\small\bmx a_1&1&0&\dots&0\\
          a_2&x_1&1&\dots&0 \\
          \vdots&\vdots && \ddots&\vdots \\
          a_{k}&x_{k-1}& \dots& &1\\
          a_{k+1}&x_{k}& \dots& \dots&x_1\\
          \emx}
          &=a_1\cdot \det T(k)
          - \det {\small\bmx a_2&1&0&\dots&0\\
          a_3&x_1&1&\dots&0 \\
          \vdots&\vdots && \ddots&\vdots \\
          a_{k}&x_{k-2}& \dots& &1\\
          a_{k+1}&x_{k-1}& \dots& \dots&x_1\\
          \emx}\\
          &=a_1m_{k} - \sum_{i=1}^{k}(-1)^{i-1}a_{i+1}m_{k-i}   \hspace{1cm} \text{(by induction)} \\
          &=a_1m_{(k+1)-1}+\sum_{i=2}^{k+1}(-1)^{i-1}a_im_{k+1-i}\\
          &=\sum_{i=1}^{k+1}(-1)^{i-1}a_im_{(k+1)-i}
\end{align*}

In order to obtain (2) one only has to solve  (1) for $x_k$ which actually is the last term in the sum up to the factor $(-1)^{k-1}$:

We obtain
 \begin{align*}
	(-1)^{k-1}x_k&=m_kx_0-\sum_{i=1}^{k-1}(-1)^{i-1}x_im_{k-i}=\sum_{i=1}^{k-1}(-1)^ix_im_{k-i}+m_kx_0\\
	x_k&=\sum_{i=0}^{k-1}(-1)^{i+k-1}x_im_{k-i}=\sum_{j=1}^k(-1)^{k-j+k-1}m_jx_{k-j},
\end{align*}
where $k-j+k-1\equiv j-1 \text{  mod } 2$.\hspace{3cm} $\Box$\\

The recursions (1) and (2) lead us directly to the {\bf involution property}. 
It will be sufficient to show that $\varphi(\varphi(x_k))=x_k$ for $1\le k\le n$. Since
\[
	 \varphi(\varphi(x_1))= \varphi(m_1))= \varphi(x_1)=m_1=x_1
\]
by induction we obtain for $1\le k\le n-1$
\begin{align*}
	\varphi(\varphi(x_{k+1}))&= \varphi(m_{k+1})= \varphi(\sum_{i=1}^{k+1}(-1)^{i-1}x_im_{(k+1)-i})\\
		&=\sum_{i=1}^{k+1} (-1)^{i-1}\varphi(x_i) \varphi(m_{(k+1)-i})\\
		&=\sum_{i=1}^{k+1}(-1)^{i-1} m_i \varphi( \varphi(x_{(k+1)-i}))\\ 
		&=\sum_{i=1}^{k+1}(-1)^{i-1}m_ix_{(k+1)-i}\hspace{ 1cm}  \text{(by induction)}\\
		&=x_{k+1}\hspace{2cm} \Box
\end{align*}

\end{document}